\documentclass[12pt,leqno]{article}
\usepackage{amsmath, amsthm, amsfonts, amssymb, color}
\usepackage{mathrsfs}
\theoremstyle{definition}

\newcommand{\scr}[1]{\mathscr #1}
\definecolor{wco}{rgb}{0.5,0.2,0.3}

\numberwithin{equation}{section} \theoremstyle{remark}

\newcommand{\ua}{\uparrow}

\title{{\bf Log-Harnack Inequality for Stochastic Burgers Equations and Applications}\footnote{Supported in
 part by WIMCS and NNSFC(10721091).}
}
\author{
{\bf Feng-Yu Wang$^{a),b)}$, Jiang-Lun Wu$^{b)}$, Lihu Xu$^{c)}$}\\
\footnotesize{$^{a)}$ School of Math. Sci. and Lab. Math. Com. Sys.,
Beijing Normal
University, Beijing 100875, China}\\
 \footnotesize{$^{b)}$ Department of Mathematics,
Swansea University, Singleton Park, SA2 8PP, UK}\\
\footnotesize{$^{c)}$ PO Box 513, EURANDOM, 5600 MB  Eindhoven, The Netherlands} \\
}
\begin{document}
\def\R{\mathbb R}  \def\ff{\frac} \def\ss{\sqrt} \def\B{\mathbf
B}
\def\N{\mathbb N} \def\kk{\kappa} \def\m{{\bf m}}
\def\dd{\delta} \def\DD{\Delta} \def\vv{\varepsilon} \def\rr{\rho}
\def\<{\langle} \def\>{\rangle} \def\GG{\Gamma} \def\gg{\gamma}
  \def\nn{\nabla} \def\pp{\partial} \def\EE{\scr E}
\def\d{\text{\rm{d}}} \def\bb{\beta} \def\aa{\alpha} \def\D{\scr D}
  \def\si{\sigma} \def\ess{\text{\rm{ess}}}
\def\beg{\begin} \def\beq{\begin{equation}}  \def\F{\scr F}
\def\Ric{\text{\rm{Ric}}} \def\Hess{\text{\rm{Hess}}}
\def\e{\text{\rm{e}}} \def\ua{\underline a} \def\OO{\Omega}  \def\oo{\omega}
 \def\tt{\tilde} \def\Ric{\text{\rm{Ric}}}
\def\cut{\text{\rm{cut}}} \def\P{\mathbb P} \def\ifn{I_n(f^{\bigotimes n})}
\def\C{\mathbb C}      \def\aaa{\mathbf{r}}     \def\r{r}
\def\gap{\text{\rm{gap}}} \def\prr{\pi_{{\bf m},\varrho}}  \def\r{\mathbf r}
\def\Z{\mathbb Z} \def\vrr{\varrho} \def\ll{\lambda}
\def\L{\scr L}\def\Tt{\tt} \def\TT{\tt}\def\II{\mathbb I}
\def\i{{\rm in}}\def\Sect{{\rm Sect}}\def\E{\mathbb E} \def\H{\mathbb H}
\def\M{\scr M}\def\Q{\mathbb Q} \def\texto{\text{o}} \def\LL{\Lambda}
\def\Rank{{\rm Rank}} \def\B{\scr B}
\def\T{\mathbb T}\def\i{{\rm i}} \def\ZZ{\hat\Z}\def\c{{\mathbf c}}

\maketitle
\begin{abstract}  By  proving an $L^2$-gradient estimate for the corresponding Galerkin approximations, the log-Harnack inequality is established for the semigroup associated to a class of stochastic Burgers equations. As applications, we derive the strong Feller property of the semigroup, the irreducibility of the solution, the entropy-cost inequality for the adjoint semigroup, and entropy upper bounds of the transition density.  \end{abstract} \noindent

 \vskip 2cm

\section{Introduction}

Let $\T= \R/(2\pi\Z)$ be equipped with the usual Riemannian metric, and let $\d\theta$ denote the Lebesgue   measure on $\T$. Then

$$\H:=\bigg\{x\in L^2(\d\theta): \int_\T x(\theta)\d\theta =0\bigg\}$$ is a separable real Hilbert space with  inner product and norm

$$\<x,y\>:=\int_\T x(\theta)y(\theta)\d\theta,\ \ \|x\|;=\<x,x\>^{1/2}.$$
For $x\in C^2(\T)$, the Laplacian operator $\DD$ is given by $\DD x= x''.$ Let $(A,\D(A))$ be the closure of $(-\DD, C^2(\T)\cap \H)$ in $\H$,
which is a positively definite self-adjoint operator on $\H$. Then

$$V:= \D(A^{1/2}),\ \<x,y\>_V:= \<A^{1/2}x, A^{1/2} y\>$$ gives rise to a Hilbert space, which is   densely and compactly embedded in $\H$.
By the integration by parts formula, for any $x\in C^2(\T)$ we have

$$\|x\|_V^2 :=\|A^{1/2}x\|^2==\int_\T(xAx)(\theta)\d\theta = \int_\T |x'(\theta)|^2\d \theta.$$
Moreover, for $x,y\in C^1(\T)\cap\H$, set $B(x,y):= xy'.$ Then $B$ extends to a unique bounded bilinear operator
$B: V\times V\to \H$ with (see Proposition
\ref{P2.1} below)

\beq\label{1.1} \|B\|_{V\to \H}:= \sup_{\|x\|_V,\|y\|_V\le 1} \|B(x,y)\|\le \ss\pi.\end{equation}
Consider the following stochastic Burgers equation

\beq\label{1.2} \d X_t= -\big\{\nu AX_t+B(X_t)\}\d t+Q\d
W_t,\end{equation} where $\nu>0$ is a constant, $B(x):=B(x,x)$ for
$x\in V$, $Q$ is a Hilbert-Schmidt operator on $\H$, and $W_t$ is
the cylindrical Brownian motion on $\H$. According to to
\cite[Chapter 5]{DP} (see also \cite[Theorem 14.2.4]{DZ}), for any
$x\in \H$, this equation has a unique solution with the initial $X_0=x$, which
is a continuous   Markov process on $\H$ and is denoted by $X_t^x$
from now on. If moreover $x\in V$, then $X_t^x$ is a continuous
process on $V$ (see Proposition \ref{P2.3} below). We are concerned with the
associated Markov semigroup  $P_t$ given by

$$P_t f(x):=\E f(X_t^x),\ \ x\in \H, t\ge 0$$ for $f\in \B_b(\H)$, the set of all bounded measurable functions on $\H$.

The purpose of this paper is to investigate regularity properties of $P_t$, such as strong Feller property, heat kernel upper bounds, contractivity properties,
 and entropy-cost inequalities. To do this, a powerful tool is the dimension-free Harnack inequality introduced in \cite{W97} for diffuions on Riemannian manifolds (see also \cite{ATW06, ATW09} for further development). In recent years, this inequality has been established and applied intensively in the study of SPDEs (see e.g. \cite{RW03, W07, LW08, DRW09, ES, W10Ann} and references within). In general, this type of Harnack inequality can be formulated as

 \beq\label{H} (P_tf)^\aa(x)\le (P_t f^\aa)(y) \exp[C_\aa(t,x,y)],\ \ f\ge 0,\end{equation} where $\aa>1$ is a constant, $C_\aa$ is a positive function on $(0,\infty)\times\H^2$ with $C_\aa(t, x,x)=0$, which is determined by the underlying stochastic equation.

 On the other hand, in some cases this kind of Harnack inequality is not available, so that the following weaker version
 (i.e. the log-Harnack inequality)

 \beq\label{LH0}   P_t\log f(x)\le \log P_t f(y)+ C(t,x,y), \ \ f\ge 1\end{equation}
 becomes an alternative tool in the study. In general, according to \cite[Section 2]{W10}, (\ref{LH0}) is the limit version of (\ref{H})
 as $\aa\to\infty.$
This inequality has been established in \cite{RW10} and \cite{W10}, respectively,  for   semi-linear SPDEs with multiplicative noise and the Neumann semigroup on non-convex manifolds.

As for the stochastic Burgers equation (\ref{1.2}), by using $A^{1+\si}$ for $\si>\ff 1 2$ to replace $A$ (i.e. the hyperdissipative equation is concerned), the first   and the third named authors   established an explicit  Harnack inequality of type (\ref{H}) in \cite{WX}, where a more general framework, which  includes also the stochastic hyperdissipative Navier-Stokes equations, was   considered. But, when $\si\le \ff 1 2$, the known arguments (i.e. the coupling argument and gradient estimate) to prove (\ref{H}) are no longer valid. Therefore, in this paper we turn to investigate the log-Harnack inequality for $P_t$ associated to  (\ref{1.2}), which also provides some important regularity properties of the semigroup (see Corollary \ref{C1.2}  below). Note that the stochastic Burgers equation does not satisfy the  Lipschitz and monotone  conditions required in \cite{RW10}, the present study can not be covered there.

To state our main result, we introduce the intrinsic norm   induced by the diffusion part of the solution. For any $x\in H$, let

$$\|x\|_Q:= \inf\big\{\|z\|_\H: z\in\H, Q^*z=x\big\},$$ where $Q^*$ is the adjoint operator of $Q$, and we take $\|x\|_Q=\infty$ if the set in the right hand side is empty.  Moreover, let $\|\cdot\|$ and $\|\cdot\|_{HS}$ denote the operator norm and the Hilbert-Schmidt norm respectively for bounded linear operators on $\H$.

\beg{thm}\label{T1.1} Assume that  $\nu^3\ge 4\pi \|A^{-1/2}Q\|^2$.
Then for any $f\in \B_b(\H)$ with $f\ge 1$,

\beq\label{LH} P_t\log f(x)\le \log P_t f(y) +\ff{2\pi
\|Q\|_{HS}^2\|x-y\|_Q^2}{1-\exp[-\ff{4\pi}{\nu^2}\|Q\|_{HS}^2
t]}\exp\Big[\ff{4\pi}{\nu^2}(\|x\|^2\lor\|y\|^2)\Big]\end{equation}
holds for $t>0$ and $x,y\in \H$. \end{thm}

Before introducing consequences of Theorem \ref{T1.1}, let us recall that the invariant probability measure of $P_t$ exists, and any invariant probability measure $\mu$ satisfies $\mu(V)=1.$ These follow immediately since $V$ is compactly embedded in $\H$ and due to the It\^o formula one has

$$\E \|X_t^0\|_H^2+ 2 \nu \int_0^t\E \|X_s^0\|_V^2\d s \le  \|Q\|_{HS}^2t,\ \ t\ge 0.$$
Next, for any two probability measures $\mu_1,\mu_2$ on $\H$, let $W_\c(\mu_1,\mu_2)$ be the  transportation-cost between them  with cost function
$$(x,y)\mapsto \c(x,y):=\|x-y\|_Q^2 \exp\Big[\ff{4\pi}{\nu^2}(\|x\|^2\lor\|y\|^2)\Big].$$ That is,

$$W_c(\mu_1,\mu_2)= \inf_{\mu\in \scr C(\mu_1,\mu_2)}\int_{\H\times\H} \c(x,y)\mu(\d x,\d y),$$ where $\scr C(\mu_1,\mu_2)$ is the
set of all couplings of $\mu_1$ and $\mu_2$. Finally, let

$$B_V(x,r)=\big\{z\in V: \|z-x\|_V<r\big\},\ \ \ x\in V, r>0.$$

\beg{cor}\label{C1.2}  Assume that  $\nu^3\ge 4\pi \|A^{-1/2}Q\|^2$.
\beg{enumerate} \item[$(1)$] For any $t>0,\ P_t$ is intrinsic strong
Feller, i.e.
$$\lim_{\|x-y\|_Q\to 0} P_t f(y)= P_t f(x),\ \ x\in \H, f\in\B_b(\H).$$ \item[$(2)$] Let $\mu$
be an invariant probability measure of $P_t$ and let $P_t^*$ be the adjoint operator of $P_t$ w.r.t. $\mu$. Then the entropy-cost inequality
$$\mu((P_t^*f)\log P_t^* f)\le \ff{2\pi \|Q\|_{HS}^2}{1-\exp[-\ff{4\pi}{\nu^2}\|Q\|_{HS}^2 t]}W_\c(f\mu,\mu),\ \ f\ge 0,\mu(f)=1$$ holds for all $t>0$.
\item[$(3)$] Let $\|\cdot\|_Q\le C\|\cdot\|_V$ hold for some constant $C>0$. Then
\beq\label{IR} \P(X_t^y\in B_V(x,r))>0,\ \ \ x,y\in V,
t,r>0.\end{equation}   Consequently, $P_t$ has a unique invariant
probability measure $\mu$, which is fully supported on $V$, i.e.
$\mu(V)=1$ and $\mu(G)>0$ for any non-empty open set $G\subset V$.
Furthermore, $\mu$ is strong mixing, i.e. for any $f \in \B_b(\H)$,
$$\lim_{t \rightarrow \infty} P_tf(x)=\mu(f), \ \ \ \forall \ x \in V.$$
\item[$(4)$] Under the same assumption as in $(3), \ P_t$ has a transition density $p_t(x,y)$ w.r.t. $\mu$ on $V$
such that the entropy inequalities \beq \label{E1}\int_V
p_t(x,z)\log \ff{ p_t(x,z)}{p_t(y,z)} \mu(\d z)\le \ff{2\pi
\|Q\|_{HS}^2\c(x,y)}{1-\exp[-\ff{4\pi}{\nu^2}\|Q\|_{HS}^2 t]}
\end{equation} and
 \beq\label{E2}\int_V p_t(x,y)\log p_t(x,y)\mu(\d y)\le -\log \int_V \exp\bigg[-\ff{2\pi \|Q\|_{HS}^2\c(x,y)}{1-\exp[-\ff{4\pi}{\nu^2}\|Q\|_{HS}^2 t]}\bigg]\mu(\d y) \end{equation} hold for all $t>0$ and $x,y\in V$.
\end{enumerate}\end{cor}

To prove the above results,  we  present in Section 2 some preparations including   a brief proof of (\ref{1.1}), a convergence theorem for the Galerkin approximation of (\ref{1.2}), and the continuity of the solution in $V$. Finally, complete proofs of Theorem \ref{T1.1} and Corollary \ref{C1.2} are addressed in Section 3.

\section{Some preparations}

 Obviously, (\ref{1.1}) is equivalent to the following result.

\beg{prp}\label{P2.1}  $\|B(x,y)\|^2\le \pi \|x\|_V^2\|y\|_V^2 $ holds for any $x,y\in V$. \end{prp}

\beg{proof} We shall take the continuous version for an element in $V$. Since $\int_\T x(\theta)\d\theta=0$, there exists $\theta_0\in \T$ such that $x(\theta_0)=0$. For any $\theta\in \T$, let $\gg : [0, d(\theta_0,\theta)]\to\T$ be the minimal geodesic from $\theta_0$ to $\theta$, where $d(\theta_0,\theta)(\le\pi)$ is the Riemannian distance between these two points. By the Schwartz inequality we have

$$|x(\theta)|^2 =\bigg|\int_0^{d(\theta_0,\theta)} \ff{\d}{\d s} x(\gg_s)\d s\bigg|^2 \le d(\theta_0,\theta) \int_\T |x'(\xi)|^2\d\xi\le\pi\|x\|_V^2.$$ Therefore,

$$\|B(x,y)\|^2=\int_\T |(xy')(\theta)|^2\d\theta\le \pi\|x\|_V^2\|y\|_V^2.$$\end{proof}

 \paragraph{Remark.} From the proof we see that (\ref{1.1}) is a   property in one-dimension, since for $d\ge 2$ there is no any constant $C\in (0,\infty)$ such that

 $$\|x\|_\infty^2\le C \int_{\T^d}|\nn x|^2(\theta)\d\theta,\ \ x\in C^1(\T^d)$$ holds.
 The invalidity of (\ref{1.1}) in high dimensions  is the main reason
why we only consider here the stochastic Burgers equation rather
than the stochastic Navier-Stokes equation.

\

 Next, due to the fact that to prove the log-Harnack inequality we have to apply the  It\^o formula for a reasonable class of reference functions which is, however, not available in infinite-dimensions,   we need to make use of the finite-dimensional approximations. To introduce the Galerkin approximation of (\ref{1.2}), let us formulate $\H$ by using the standard ONB
 $\{e_k: k\in \Z\}$ for the complex Hilbert space $L^2(\T\to\C; \d\theta),$ where

 $$e_k(\theta):= \ff 1 {\ss{2\pi}} \e^{\i k\theta},\ \ \ \theta\in \T.$$ Obviously, $\DD e_k=-k^2 e_k$ holds for all $k\in\Z$, and an element

 $$x:=\sum_{k\in \Z} x_k e_k,\ \ x_k\in\C$$ belongs to $\H$ if and only if $x_0=0$, $\bar x_k=x_{-k}$   for
 all $k\in \hat\Z:= \Z\setminus\{0\},$  and $\sum_{k\in\hat \Z} |x_k|^2<\infty.$ Therefore,

 $$\H= \bigg\{\sum_{k\in\hat\Z} x_ke_k: \bar x_k=x_{-k}, \sum_{k\in\hat \Z} |x_k|^2<\infty\bigg\}.$$For any $m\in\N$, let

 $$\H_m=\big\{x\in \H:\ \<x,e_k\>=0\ \text{for}\ |k|>m\big\},$$ which is a finite-dimensional Euclidean space.
 Let $\pi_m: \H\to\H_m$ be the orthogonal projection. Let $B_m=\pi_m B$ and $Q_m=\pi_m Q$. Consider the following  stochastic differential
 equation on $\H_m$:

 \beq\label{2.1} \d X_t^{(m)}= -\big\{\nu AX_t^{(m)} +B_m(X_t^{(m)})\big\}\d t+Q_m\d W_t.\end{equation} Since coefficients in this equation are smooth and

 $$\d\|X_t^{(m)}\|^2\le 2\|Q_m\|_{HS}^2\d t+ 2\<X_t^{(m)}, Q_m\d W_t\>,$$ we conclude that starting from  any $x\in \H_m$ this equation has a unique strong solution $X_t^{m,x}$  which is non-explosive. Let

 $$P_t^{(m)} f(x)=\E f(X_t^{m,x}),\ \ t\ge 0, x\in\H_m, f\in\B_b(\H_m).$$
 In the spirit of   \cite[Theorem 5.7]{DP},  the next result implies

 \beq\label{2.2} P_tf(x)= \lim_{m\to\infty} P_t^{(m)} f(  x_m),\ \ x\in\H, f\in C_b(\H) \end{equation}   for $\{x_m\in\H_m\}_{m\ge 1}$ such that
 $x_m\to x$ in $\H$.

 \beg{prp}\label{P2.2}   For any $\{x_m\in\H_m\}_{m\ge 1}$ such that $\|x-x_m\|_H\to 0$, we have $\|X_t^x-X_t^{m,x_m}\|\to 0$ in probability as $m\to\infty$. Consequently, $(\ref{2.2})$ holds.\end{prp}
\beg{proof}  Simply denote $X_t(m)=X_t^{m,x_m}$ and $X_t=X_t^x$.   It is easy to see that

\beq\label{NN}\E\int_0^t (\|X_s\|_V^2+\|X_s(m)\|_V^2)\d s \le C(1+t)\end{equation} holds for some constant $C>0.$ By the It\^o formula we have

\beq\label{A1}\beg{split}&\|X_t-X_t(m)\|^2\\
&= -2 \int_0^t\big\{\nu \|X_s-X_s(m)\|_V^2+ \<B(X_s) -B(X_s(m)),
X_s-X_s(m)\>\big\}\d s +\eta_t(m),
\end{split}\end{equation}
where
$$\eta_t(m):=\|Q-Q_m\|_{HS}^2t + \|x-x_m\|^2 +2\sup_{r\in [0,t]}\bigg|\int_0^r\<X_s-X_s(m), (Q-Q_m)\d W_s\>\bigg|,$$ which goes to $0$ as $m\to\infty$. Since by (\ref{1.1})

\beg{equation*}\beg{split} &|\<B(x)-B(y),x-y\>|=|\<B(x,x-y)+B(x-y,y),x-y\>|\\
&\le \pi
\|x-y\|(\|x\|_V+\|y\|_V)\|x-y\|_V,\end{split}\end{equation*} it
follows from (\ref{A1}) that

$$ \|X_r-X_r(m)\|^2\le \ff {\pi}{\nu}\int_0^r \|X_s-X_s(m)\|^2(\|X_s\|_V^2+\|X_s(m)\|_V^2)\d s +\eta_t(m),\ \ r\in [0,t].$$
Therefore,

$$\|X_t-X_t(m)\|^2\le \eta_t(m) \exp\bigg[\ff {\pi}{\nu}\int_0^t  (\|X_s\|_V^2+\|X_s(m)\|_V^2)\d s\bigg].$$  Combining this with (\ref{NN}) we obtain that for any $N>0$,

$$\P\big(\|X_t-X_t(m)\|^2\ge \eta_t(m) \e^{N\pi/\nu}\big)\le \P\bigg(\int_0^t  (\|X_s\|_V^2+\|X_s(m)\|_V^2)\d s\ge N\bigg)\le \ff{C(1+t)}N$$ which goes to $0$ as $N\to\infty$. Since $\eta_t(m)\to 0$ as $m\to\infty$, this implies that $\|X_t-X_t(m)\|\to 0$ in probability as $m\to\infty.$
\end{proof}

Finally, we have the following  result for the continuity of the solution in $V$.

\beg{prp}\label{P2.3} For any $x\in V$,  $X_t^x$ is a continuous process in $V$.  \end{prp}

\beg{proof} For fixed  $x\in V$ and $T>0$,  we introduce the map

$$Y: C([0,T];V)\to C([0,T];V),$$
such that for any $u\in C([0,T];V)$, $\{Y_t(u)\}_{t\in [0,T]}$ solves the deterministic equation

\beq\label{2.3} \dot Y_t(u)= -\big\{\nu AY_t(u)
+B(Y_t(u)+u_t)\big\},\ \  Y_0(u)=x.\end{equation} Then  $Y(u)\in
C([0,T];V)$,  see e.g. \cite [Theorem 3.2]{Te} (the theorem is  for
2D Navier-Stokes equation,  and the proof works also  for our case).

Next, let

$$Z_t=\int_0^t \e^{-\nu(t-s)A}Q\d W_s.$$Since $Q$ is Hilbert-Schmidt  on $\H$, $Z_t$  is a continuous process on $V$ (see e.g. \cite [Theorem 5.9]{DZ92}).  Therefore,
$ X_t^x= Y_t(Z)+Z_t$ is also continuous in $V$.\end{proof}

\section{Proofs of Theorem \ref{T1.1} and Corollary \ref{C1.2}}

According to \cite{RW10}, the key step to prove the log-Harnack inequality for $P_t^{(m)}$ is the $L^2$-gradient estimate 

$$|Q_m DP_t^{(m)} f|^2(x)\le (P_t^{(m)}|Q_m Df|^2)(x) C(t,x),\ \ f\in C_b^1(\H_m)$$ for some continuous function $C$ on $(0,\infty)\times \H_m,$  where $D$ is the gradient operator on $\H_m$, i.e. for any $f\in C^1(\H_m),$ the element $Df(x)\in\H_m$ is determined by

$$\<Df(x),h\>=D_h f(x):=\lim_{\vv\to 0} \ff{f(x+\vv h)-f(x)}{\vv},\ \ h\in \H_m.$$
To derive the desired gradient estimate, we need the following lemma.

\beg{lem}\label{L3.1} For any $x\in\H_m$ and $t\ge 0$,

$$\E\exp\bigg[\ff {\nu}{2\|A^{-1/2}Q_m\|^2}\bigg(\|X_t^{m,x}\|^2 +\nu\int_0^t \|X_s^{m,x}\|_V^2\d s\bigg)\bigg]\le
\exp\bigg[\ff {\nu(\|x\|^2+\|Q_m\|_{HS}^2
t)}{2\|A^{-1/2}Q_m\|^2}\bigg].$$\end{lem}

\beg{proof} By the It\^o formula, we have \beq\label{3.1}
\d\|X_t^{m,x}\|^2 + 2\nu\|X_t^{m,x}\|_V^2\d t = \|Q_m\|_{HS}^2\d t +
2\<X_t^{m,x}, Q_m \d W_t\>.\end{equation} Let

$$\tau_n=\inf\big\{t\ge 0:  \|X_t^{m,x}\|\ge n\},\ \ n\in\N.$$ We have $\tau_n\to\infty$ as $n\to\infty$. Let

$$M_t^{(n)}= \int_0^{t\land \tau_n} \<X_s^{m,x},Q_m\d W_s\>.$$ Then for any   $\ll>0$

$$t\mapsto\exp\big[2\ll M_t^{(n)}-2\ll^2 \<M^{(n)}\>_t\big]$$ is a martingale. Therefore,   it follows from (\ref{3.1}) that

\beq\label{3.2} \beg{split} &\E
\exp\bigg[\ll\|X_{t\land\tau_n}^{m,x}\|^2 + 2\nu\ll
\int_0^{t\land\tau_n}\|X_s^{m,x}\|_V^2\d s- 2\ll^2
\int_0^{t\land\tau_n}
\|Q_m^* X_s^{m,x}\|^2\d s\bigg]\\
&\le \E\exp\bigg[\ll\big(\|x\|^2 +t\|Q_m\|_{HS}^2\big) + 2\ll M_t^{(n)} -2 \ll^2\<M^{(n)}\>_t\bigg]\\
&= \exp\big[\ll(\|x\|^2+t\|Q_m\|_{HS}^2)\big].\end{split}\end{equation}
Noting that

$$\|Q_m^*x\|=\|Q_m^*A^{-1/2}A^{1/2}x\|\le \|Q_m^*A^{-1/2}\|\cdot \|x\|_V =\|A^{-1/2}Q_m\|\cdot\|x\|_V,\ \ x\in\H_m,$$ by letting $n\uparrow\infty$ in (\ref{3.2}) and taking

$$\ll=\ff \nu{2\|A^{-1/2}Q_m\|^2},$$ we complete the proof.
\end{proof}

\beg{lem}\label{L3.2} Let $\nu^3\ge 4\pi \|A^{-1/2}Q_m\|^2.$ Then
for any $f\in C_b^1(\H_m),$

$$\|Q_m DP_t^{(m)} f\|^2(x)\le \big(P_t^{(m)}\|Q_m Df\|^2\big)(x) \exp\bigg[\ff{2\pi}{\nu^2} \big(\|x\|^2+t\|Q_m\|_{HS}^2\big)\bigg]$$ holds for all $t\ge 0$ and $x\in\H_m.$\end{lem}

\beg{proof} Let $h\in \H_m$. According to e.g. \cite[Section 5.4]{DP},

$$D_hX_t^{m,x}:= \lim_{\vv\to 0}\ff{X_t^{m,x+\vv h}-X_t^{m,x}}\vv,\ \ \ t\ge 0$$ exists and solves the ordinary differential equation

$$\ff{\d}{\d t} D_h X_t^{m,x}=-\big\{\nu AD_h X_t^{m,x} +\tt B_m (X_t^{m,x}, D_h X_t^{m,x})\big\},$$ where $\tt B_m (x,y):= B(x,y)+B(y,x)$ for $x,y\in\H_m.$
By (\ref{1.1}),  this implies that

\beg{equation*}\beg{split} \ff{\d}{\d t} \|D_h X_t^{m,x}\|_V^2 &= -2\nu \|AD_h X_t^{m,x}\|^2 - 2 \<A D_h X_t^{m,x}, \tt B_m (X_t^{m,x}, D_h X_t^{m,x})\>\\
&\le \ff 1 {2\nu} \|\tt B_m(X_t^{m,x}, D_hX_t^{m,x})\|^2\le
\ff{2\pi}{\nu}\|X_t^{m,x}\|_V^2\|D_hX_t^{m,x}\|_V^2.\end{split}\end{equation*}
Therefore,

$$\|D_h X_t^{m,x}\|_V^2 \le \|h\|_V^2\exp\bigg[\ff{2\pi}{\nu}\int_0^t\|X_s^{m,x}\|_V^2\d s\bigg].$$ Since $\nu^3\ge 4\pi \|A^{-1/2}Q_m\|^2$ implies that

$$\ff{\nu^2}{2\|A^{-1/2}Q_m\|^2}\ge \ff{2\pi}\nu,$$ by Lemma \ref{L3.1} and using the Jensen inequality we arrive at

\beq\label{3.3} \E\|D_h X_t^{m,x}\|_V^2\le \|h\|_V^2
\exp\bigg[\ff{2\pi}\nu
\big(\|x\|^2+t\|Q_m\|_{HS}^2\big)\bigg].\end{equation} Consequently,
by the dominated convergence theorem we obtain

\beq\label{3.4} D_h P_t^{(m)} f (x)= \E\<Df(X_t^{m,x}), D_h X_t^{m,x}\>,\ \ f\in C_b^1(\H_m), x\in \H_m, t\ge 0.\end{equation}
On the other hand, we have

\beq\label{3.5} \beg{split}\|Q_m D P_t^{(m)}f\|^2 &= \sup_{\|\tt h\|\le 1} \<Q_m DP_t^{(m)} f, \tt h\>^2 = \sup_{\|\tt h\|\le 1} \<DP_t^{(m)}f, Q_m^* \tt h\>^2\\
&=\sup_{\|h\|_{Q_m}\le 1} |D_hP_t^{(m)}f|^2,\end{split}\end{equation} where

$$\|h\|_{Q_m}:= \inf\{\|z\|: z\in\H_m, Q_m^* z=h\}$$ and $\|h\|_{Q_m}=\infty$ if the set on the right hand side is empty. Now, for any $h\in \H_m$ with
 $\|h\|_{Q_m}\le 1$, let $\{z_n\}_{n\ge 1}\subset \H$ be such that $Q_m^*z_n=h$ and $\|z_n\|\le 1+\ff 1 n.$ By (\ref{3.4}) we have

 \beg{equation*}\beg{split} &|D_hP_t^{(m)}f|^2(x)=\big(\E\<Df(X_t^{m,x}), D_h X_t^{m,x}\>\big)^2 =\big(\E\<Q_mDf(X_t^{m,x}), D_{z_n} X_t^{m,x}\>\big)^2\\
 &\le \big(\E\|Q_m Df(X_t^{m,x})\|^2\big)\E\|D_{z_n} X_t^{m,x}\|^2= \big(\E\|Q_m Df(X_t^{m,x})\|^2\big)\E\|D_{A^{-1/2} z_n} X_t^{m,x}\|_V^2.\end{split}\end{equation*}
 Combining this with (\ref{3.3}) and  (\ref{3.5}) and  letting $n\uparrow\infty$, we complete the proof.    \end{proof}

According to the $L^2$-gradient estimate in Lemma \ref{L3.2}, we are able to prove the log-Harnack inequality for $P_t^{(m)}$ as in \cite{RW10}.

\beg{prp}\label{P3.3} Let $\nu^3\ge 4\pi \|A^{-1/2}Q_m\|^2.$ For any
$f\in \B_b(\H_m)$ with $f\ge 1$,

$$P_t^{(m)}\log f(x)\le \log P_t^{(m)} f(y) +
\ff{\pi\|Q_m\|_{HS}^2\|x-y\|_{Q_m}^2\exp[\ff{2\pi}{\nu^2}(\|x\|^2\lor\|y\|^2)]}{1-\exp[-\ff{2\pi}{\nu^2}\|Q_m\|_{HS}^2t]}$$
holds for all $t>0$ and $x,y\in \H_m.$\end{prp}

\beg{proof} It suffices to prove for $\|x-y\|_{Q_m}<\infty$. Let $\{z_n\}\subset\H_m$ be such that
$Q_m^* z_n=x-y$ and $\|z_n\|^2\le \|x-y\|_{Q_m}^2 +\ff 1 n.$  Let $\gg\in C^1([0,t];\R)$ such that $\gg(0)=0,\gg(t)=1.$ Finally, let $x_s= (x-y)\gg(s)+y,\ s\in [0,t].$ Then, by Lemma \ref{L3.2} we have (see \cite[Proof of Theorem 2.1]{RW10} for explanation of the second equality)

\beg{equation*}\beg{split}& P_t^{(m)}\log f(x)-\log P_t^{(m)} f(y)  = \int_0^t \ff{\d}{\d s} \big\{P_s^{(m)}\log P_{t-s}^{(m)} f\big\}(x_s)\d s\\
&=\int_0^t\Big\{-\ff 1 2P_s^{(m)} \|Q_m D \log P_{t-s}^{(m)} f\|^2 +\gg'(s) \<x-y, DP_s^{(m)}\log P_{t-s}^{(m)}f\>\Big\}(x_s)\d s\\
&\le \int_0^tP_s^{(m)} \Big\{-\ff 1 2 \|Q_m D \log P_{t-s}^{(m)} f\|^2\\
&\qquad\qquad\qquad +|\gg'(s)| \cdot\|z_n\|\e^{2\pi(\|x_s\|^2
 +\|Q_m\|_{HS}^2s)/\nu^2}\|Q_m D\log P_{t-s}^{(m)}f\|\Big\}(x_s)\d s\\
&\le \ff{\|z_n\|^2} 2 \int_0^t |\gg'(s)|^2
\e^{4\pi(\|x_s\|^2+\|Q_m\|_{HS}^2s)/\nu^2}\d
s.\end{split}\end{equation*} Since $\|x_s\|\le \|x\|\lor\|y\|$, by
taking

$$\gg(s)= \ff{ 1-\exp[-\ff{4\pi}{\nu^2} \|Q_m\|_{HS}^2s]}{1-\exp[-\ff{4\pi}{\nu^2} \|Q_m\|_{HS}^2 t]},\ \ s\in [0,t]$$
we obtain
$$P_t^{(m)}\log f(x)-\log P_t^{(m)} f(y)+\ff{2\pi \|Q\|_{HS}^2\|z_n\|^2}{1-\exp[-\ff{4\pi}{\nu^2}\|Q\|_{HS}^2 t]}\exp\Big[\ff{4\pi}{\nu^2}(\|x\|^2\lor\|y\|^2)\Big].$$
This completes the proof by letting $n\to\infty$.\end{proof}

\beg{proof}[Proof of Theorem \ref{T1.1}]  It suffices to prove for
$f\in C_b(\H)$ with $f\ge 1$. Let $\|x-y\|_Q<\infty$. For any
$\vv>0$, let $z\in \H$ such that $Q^*z=x-y$ and $\|z\|^2\le
\|x-y\|_Q^2+\vv.$ For any $m\in\N$, we have $Q^*_mz=\pi_mx-\pi_my$.
Let $x_m=\pi_m x, z_m=\pi_m z$ and $y_m= \pi_m y +Q^*_m(z-\pi_m z)$.
Then $z_m\in \H_m$ and $Q^*_mz_m=x_m-y_m$, so that

$$\|x_m-y_m\|^2_{Q_m}\le \|z_m\|^2\le \|x-y\|_Q^2+\vv.$$ Moreover, it is easy to see that $x_m\to x$ and $y_m\to y$ hold in $\H$. Combining these with Proposition \ref{P3.3} and (\ref{2.2}), and letting $m\to\infty$ and $\vv\to 0$, we complete the proof.
\end{proof}

\beg{proof}[Proof of Corollary \ref{C1.2}] The intrinsic strong Feller property follows from
 \cite[Proposition 2.3]{W10}, while the entropy-cost inequality in (2) follows from the proof of
Corollary 1.2 in \cite{RW10}. So, it remains to prove (3) and (4).

(a)  Applying (\ref{LH}) to $f:= 1+m1_{B(x,r)}$ for $m\ge 1$, we
obtain

\beq\label{WF2} P_t\log (1+m1_{B_V(x,r)})(x) \le \log\big\{1+m P_t
1_{B_V(x,r)}(y)\big\}+\aa(t)\c(x,y),\  \ t>0, m\ge 1\end{equation}
for some function $\aa: (0,\infty)\to (0,\infty)$ independent of
$x,y$ and $m$. By Proposition \ref{P2.3} we have $\|X_t^x-x\|_V\to
0$ as $t\to 0$. Then there exists $t_0>0$ depending only on $x$ such
that

$$\P(\|X_t^x-x\|_V<r)\ge \ff 1 2,\ \ \ t\in [0,t_0].$$ Thus, if $\P(X_t^y\in B_V(x,r))=0$ for some $t\in (0,t_0]$,  then (\ref{WF2}) yields that

 $$\ff 1 2 \log(1+m)\le P_t\log (1+m1_{B_V(x,r)})(x) \le \aa(t)\c(x,y),\ \ \ m\ge 1,$$ which is impossible since $\|\cdot\|_Q\le C\|x-y\|_V$ implies that
 $\c(x,y)<\infty$ for $x,y\in V$. Therefore,

 $$\P(X_t^z\in B_V(x,r))>0, \ \ \ t\in (0,t_0], z\in V.$$ Combining this with the Markov property we see that for $t>t_0$,

 $$\P(X_t^y\in B(x,r))= \int_V \P(X_{t_0}^z \in B(x,r))P_{t-t_0}(y,\d z)>0,$$ where $P_{t-t_0}(y,\d z)$ is the distribution of $X_{t-t_0}^y$.
 Therefore, (\ref{IR}) holds.

 (b)   Since (1) and $\|\cdot\|_Q\le C\|\cdot\|_V$ imply the strong Feller property of $P_t$ on $V$, by the Doob  Theorem,
 see e.g. \cite[Theorem 4.2.1]{DZ}, $P_t$ has a unique invariant measure $\mu$  on $V$. The full support
property of $\mu$, together with the strong Feller of $P_t$, implies
the existence of transition density $p_t(x,y).$  Finally, due to
\cite[Proposition 2.4(2)]{W10},
   (\ref{E1}) is equivalent to the log-Harnack inequality (\ref{LH}), while  (\ref{E2}) follows from (\ref{LH}) according to  the proof of
   \cite[Corollary 1.2]{RW10}.\end{proof}

\end{document}